\documentclass[12pt,reqno]{amsart}

\usepackage{amsmath,amsthm,amsfonts,amssymb,enumerate,color,pifont}

\usepackage{pdfpages}

\newtheorem{thm}{Theorem}[section]
\newtheorem{lemma}[thm]{Lemma}

\newtheorem{prop}[thm]{Proposition}

\theoremstyle{definition}

\newtheorem{defn}[thm]{Definition}

\newtheorem{rmk}[thm]{Remark}

\newtheorem{ack}{Acknowledgments\!\!}  

\newtheorem{Pf}{Proof$\!\!$}         
         
\newenvironment{pf}{\begin{Pf}}{\qed\end{Pf}}

\renewcommand{\qedsymbol}{{\hspace*{\fill}\rule{2.7mm}{2.7mm}}}

\newcommand{\kk}{\Bbbk}

\newcommand{\cc}{\mathbb{C}}
\newcommand{\C}{\mathbb{C}}

\newcommand{\PP}{\mathbb{P}}

\newcommand{\Z}{\mathbb{Z}}

\newcommand{\hh}{\mathcal{H}}

\newcommand{\V}{\mathcal{V}}

\newcommand{\fsl}{\mathfrak{sl}}
\newcommand{\slc}{\mathfrak{sl}^c_2}

\renewcommand{\b}{\mathfrak{b}}
\newcommand{\g}{\mathfrak{g}}
\newcommand{\h}{\mathfrak{h}}
\newcommand{\fslll}{\mathfrak{sl}(1|1)}

\newcommand{\onto}{\protect{~-\!\!\!\twoheadrightarrow}}

\begin{document}
\baselineskip21pt

\title[Associating Geometry to $\fslll$ and to $\slc(\kk)$]{Associating
Geometry to the Lie Superalgebra $\fslll$\\[2mm]
and to the Color Lie algebra $\slc(\kk)$}

\author[]{}

\subjclass[2010]{14A22, 17B70, 17B75}%

\keywords{Lie algebra, Lie superalgebra, color Lie algebra, line module}

\author[S.J.~Sierra, \v{S}.~\v{S}penko, M.~Vancliff, P.~Veerapen, E.~Wiesner]{}
\maketitle

\begin{center}
Susan J.~Sierra%
\footnote{S.J.~Sierra \& \v{S}.~\v{S}penko were supported in part by EPSRC 
grant EP/M008460/1.\\[-24pt]},
\setcounter{footnote}{0}
\v{S}pela \v{S}penko%
\footnote{}
\footnote{Some of this work was completed while \v{S}.~\v{S}penko was a
post-doctoral scholar at the University of Edinburgh.},
Michaela Vancliff%
\hspace*{0.5mm}%
\footnote{M.~Vancliff was supported in part by NSF grant DMS-1302050.},\\%
Padmini Veerapen and Emilie Wiesner
\end{center}

\vspace*{0.1in}

\begin{abstract}
\baselineskip15pt
In the 1990s, in work of Le Bruyn and Smith and in work of Le Bruyn and Van den
Bergh, it was proved that point modules and line modules over the homogenization
of the universal enveloping algebra of a finite-dimensional Lie algebra describe 
useful data associated to the Lie algebra (\cite{LS,LV}).  In particular, in the
case of the Lie algebra $\fsl_2(\C)$, there is a correspondence between Verma 
modules and certain line modules that associates a pair $(\h,\,\phi)$, where $\h$
is a two-dimensional Lie subalgebra of $\fsl_2(\C)$ and $\phi \in \h^*$
satisfies $\phi([\h, \, \h]) = 0$, to a 
particular type of line module.  In this article, we prove analogous results for the 
Lie superalgebra $\fslll$ and for a color Lie algebra associated to the Lie 
algebra $\fsl_2$.  
%
\end{abstract} 

\baselineskip18pt

\bigskip
\bigskip


\section*{Introduction}
The ideas of algebraic geometry that were introduced in \cite{ATV1} often provide 
useful tools in studying noncommutative algebras, particularly Artin-Schelter 
regular algebras. In this setting, one can associate geometric objects (points, 
lines, etc.) of an appropriate projective space to certain classes of graded modules 
(namely, point modules, line modules, etc.) over the regular algebra.

Le Bruyn and Smith in \cite{LS} and, later,  Le Bruyn and Van den Bergh in \cite{LV}
proved that this framework can be applied to Lie algebras.  In 
\cite{LS}, Le Bruyn and Smith studied the homogenization $\hh (\fsl_2(\cc))$ of the 
universal enveloping algebra of the Lie algebra $\fsl_2(\cc)$, and they showed
that most of the line modules over $\hh (\fsl_2(\cc))$ can be realized as 
homogenizations of Verma modules over $\fsl_2(\cc)$.    Taking this further, 
in~\cite{LV}, Le~Bruyn and Van den Bergh showed that these results extend to 
finite-dimensional Lie algebras~$\g$ of dimension~$n$. In particular, they showed a 
correspondence between certain $d$-dimensional linear subschemes of $\PP^n$ and 
modules induced from $d$-codimensional subalgebras of~$\g$. In the special case
of $\g = \fsl_2$, where $n = 3$ and $d = 1$, this work associates a pair 
$(\h,\,\phi)$, where $\h$ is a two-dimensional Lie subalgebra of $\g$ and $\phi 
\in \h^*$ is such that $\phi([\h, \, \h]) = 0$, to a particular type of line module 
over $\mathcal H (\fsl_2)$.  

The results in \cite{LS,LV} for Lie algebras suggest that similar techniques for 
analogous algebras -- such as Lie superalgebras and color Lie algebras -- might be 
productive.  In this article, we investigate this idea for the Lie superalgebra 
$\fslll$ and for a color Lie algebra $\slc$ associated to $\fsl_2$.  
In particular, we prove in Theorem~\ref{Lhphi}, Proposition~\ref{mainthm0}, 
Theorem~\ref{mainthm1} and Theorem~\ref{mainthm2} that, 
for $\g \in \{\fslll , \, \slc\}$, 
there is a one-to-one correspondence between the collection of certain line modules 
over an appropriate graded algebra associated to $\g$ and pairs $(\h,\,\phi)$,
where $\h$~is a two-dimensional subalgebra of $\g$ and $\phi \in \h^*$
satisfies certain conditions. 
Moreover, most line modules can be realized as homogenizations 
of modules induced from one-dimensional $\h$-modules.
We conclude the paper by discussing the challenge in applying similar ideas to 
the Lie superalgebra $\fsl(2|1)$.

\bigskip
\bigskip


\begin{ack}
This project was initiated at the first WINART (Women in Noncommutative 
Algebra and Representation Theory) workshop held at the Banff International 
Research Station (BIRS) in the Spring of 2016.  The environment at BIRS 
was particularly stimulating during the workshop and, for that, the authors 
warmly thank BIRS and the organizers and the participants of the workshop. 
\end{ack}

\bigskip
\bigskip


\section{The Lie Algebra $\fsl_2(\cc)$}

In this section, we establish notation to be used throughout the paper, and
summarize the results of \cite{LS,LV} to be generalized to the Lie superalgebra 
$\fslll$ and to the color Lie algebra $\slc$ discussed in the Introduction.

\bigskip

Throughout the article, $\kk$~denotes an algebraically closed field, and, 
unless otherwise indicated, char$(\kk) \neq 2$. 
We use $M(n,\,\kk)$ to denote the vector space of $n \times n$ 
matrices with entries in~$\kk$. For a graded $\kk$-algebra~$B$, the span of the 
homogeneous elements in~$B$ of degree~$i$ will be denoted $B_i$, and the dual of any 
vector space $V$ will be denoted $V^*$.  
For homogeneous polynomials $f_1, \ldots, f_m$, we use $\V(f_1, \ldots , f_m)$ to 
denote their zero locus in projective space.

\bigskip

In \cite{LS}, Le Bruyn and Smith show how noncommutative algebraic geometry, in the 
spirit of Artin, Tate, and Van den Bergh (cf.\ \cite{ATV1}), can be applied to the Lie 
algebra $\fsl_2(\cc)$. In particular, since the homogenization~$A$ of the universal 
enveloping algebra of~$\fsl_2(\C)$ by a central element $t$ is an Artin-Schelter 
regular algebra, the authors prove that there is a one-to-one correspondence between 
the collection of line modules over~$A$ on which $t$ acts without torsion and 
the collection of pairs $(\b, \, \lambda)$, where $\b$ is a Borel subalgebra 
of~$\fsl_2(\C)$ and $\lambda \in \C$. This result is stated in 
Theorem~\ref{thmLS} below.  To this end, consider 
the Lie algebra $\fsl_2(\C)$ with basis $\{e,\, f,\, h\}$ and Lie bracket 
\[ [e,\, f] = h, \qquad [h,\, e] = 2e, \qquad [h,\, f] = -2f.\] 
Since $\fsl_2(\C)$ embeds in  $M(2,\,\C)$ via the map 
\[
e\mapsto \begin{bmatrix} 0 & 1 \\0 & 0 \end{bmatrix}, \quad 
f\mapsto \begin{bmatrix} 0 & 0 \\1 & 0 \end{bmatrix}, \quad 
h\mapsto \begin{bmatrix} 1 & 0 \\0 & -1 \end{bmatrix},
\]
we may define the determinant function, $\det$, on elements of $\fsl_2(\C)$.
Identifying $\PP^3$ with $\PP((\kk e \oplus \kk f \oplus \kk h \oplus \kk t)^*)$,
we define the pencil of quadrics $Q(\delta) = 
\V(\det + \delta^2 t^2)$, for all $\delta \in \PP^1$
(where $Q(\infty) = \V(t^2)$), and lines 
$l_{(\b , \, \lambda)}$ in $\PP^3$ by
\[
l_{(\b , \, \lambda)} = \mathcal{V}(E, \, H - \lambda t),
\] 
where $\lambda \in \C$ and  $\mathfrak{b}$ is a Borel subalgebra of 
$\fsl_2(\C)$ with standard basis $\{E, \, H\}$. As observed in~\cite{LS},
$l_{(\b , \, \lambda)}$ is independent of the choice of standard 
basis.

In order to define a line module over~$A$, we first note that, by viewing 
$e$, $f$, $h$, $t \in A$ as having degree one, the algebra~$A$ is a 
graded $\C$-algebra.  

\begin{defn}\cite{ATV1} \  
A {\em line module} over~$A$ is a graded cyclic $A$-module with Hilbert series 
$H(x) = 1/(1-x)^2$.
\end{defn}

\begin{thm} \cite[Theorems 1 and 2]{LS}\label{thmLS}
As above, let $A$~denote the homogenization of the universal enveloping algebra 
of $\fsl_2(\C)$ by a central element $t$.\\[-7mm]
\begin{enumerate}
\item[{\rm (}a{\rm )}] 
The lines that lie on the quadrics $Q(\delta) = \V(\det + \delta^2 t^2)$ 
$($$\delta \in \PP^1${\rm )} are 
\begin{enumerate}
\item[{\rm (}i\/{\rm )}] the lines on the plane $\V(t)$, and
\item[{\rm (}ii\/{\rm )}] the lines $l_{(\b , \, \lambda)}$ where $\b$ is a 
                         Borel subalgebra and $\lambda \in \C$.
\end{enumerate}
\item[{\rm (}b{\rm )}]
The lines in $\PP^3$ that determine the line modules over $A$ are precisely the 
lines on the quadrics $Q(\delta)$, for all $\delta \in \mathbb{P}^1$.
\item[{\rm (}c{\rm )}]
The line modules over~$A$ are of two types$:$ 
\begin{enumerate}
\item[{\rm (}i\/{\rm )}] those corresponding to the lines on the plane~$\V(t)$,
that is, $A/(At + Aa)$ for all $a \in A_1 \setminus \C t$, and 
\item[{\rm (}ii\/{\rm )}] those corresponding to the lines $l_{(\b , \, \lambda)}$ 
where $\b$ is a Borel subalgebra and $\lambda \in \C$, that is,
$A/(AE + A(H-\lambda t))$, where $\{E, \, H \}$ is a standard basis of~$\b$.
\end{enumerate}
\end{enumerate}
\end{thm}

\begin{rmk} \label{borel}
For comparison with the techniques in~\cite{LV} and with our methods in 
Sections~\ref{sec:sl11} and~\ref{sec:slc}, we note that $\lambda \in \C$ defines a 
linear functional $f_{\lambda} \in \mathfrak{b}^{\ast}$ by $f_{\lambda}(E) = 0$ 
and $f_{\lambda}(H) = \lambda$, where $\b$, $E$ and~$H$ are as above.
In particular, $f_{\lambda}([\b, \, \b] ) = 0$ 
and $\lambda$ defines a 1-dimensional $\b$-module via~$f_{\lambda}$. Moreover, 
given a Borel subalgebra~$\b$ of~$\fsl_2(\C)$, any $f \in \b^*$ such that 
$f([\b, \, \b] ) = 0$ satisfies $f = f_{\lambda}$, for some $\lambda \in \C$.  
\end{rmk}

\bigskip
\bigskip


\section{The Lie Superalgebra $\fslll$} \label{sec:sl11}

In this section, we consider the Lie superalgebra $\fslll$.
We generalize to $\fslll$ the relationship between line modules over a 
certain graded algebra and subalgebras of Lie algebras given in \cite{LS,LV}.
In particular, for two-dimensional subalgebras $\h$ of $\fslll$ and linear 
functionals $\phi \in \h^{\ast}$, we show that a correspondence, analogous to 
that for $\fsl_2(\C)$, between pairs $(\h,\, \phi)$ and line modules exists.
However, only for certain pairs $(\h,\, \phi)$ and certain line modules can the 
correspondence for $\fslll$ be realized via a homogenization and 
dehomogenization process similar to that used in \cite{LS,LV}. Our main 
results are Theorem~\ref{Lhphi}, Proposition~\ref{mainthm0} and 
Theorem~\ref{mainthm1}.

In this section, we write $\g$ for $\fslll$.  Recall that $\g$ is the set of 
all $2 \times 2$ matrices with supertrace equal to zero; that is
$$
\g = \left\{ 
\begin{bmatrix} a & b \\c & d \end{bmatrix} \in M(2,\,\kk) \ : \ a-d=0 \right\}.
$$
There is a $\Z_2$-grading on $\g$ that is given by 
$\g =\g_{\overline 0} \oplus \g_{\overline 1}$ where 
$$
\g_{\overline 0}= \left\{ 
\begin{bmatrix} a & 0 \\0 & d \end{bmatrix} : a,\, d \in \kk, \ a=d
\right\}\quad \text{and}\quad
\g_{\overline 1}= \left\{ 
\begin{bmatrix} 0 & b \\c & 0 \end{bmatrix} : b,\, c \in \kk
\right\}.
$$
The Lie (super)bracket is given by 
$[X,\,Y]=XY-(-1)^{|X||Y|} YX$ for homogeneous elements 
$X,\, Y \in \g$, where $|X|$ and $|Y|$ indicate the $\Z_2$-degree of $X$ and
$Y$ respectively. 
It follows that $\g$ has basis 
$$
e=\begin{bmatrix} 0 & 1 \\0 & 0 \end{bmatrix}, \quad 
f=\begin{bmatrix} 0 & 0 \\1 & 0 \end{bmatrix}, \quad 
h=\begin{bmatrix} 1 & 0 \\0 & 1 \end{bmatrix},
$$
with relations
$$
[e,\,f]=h, \quad 0=[h,\,e]=[h,\,f]=[e,\,e]=[f,\,f]=[h,\,h].
$$
We call any subspace of~$\g$ that is closed under the Lie (super)bracket 
a {\em subalgebra} of~$\g$.


The universal enveloping algebra, $U(\g)$, of $\g$ is defined to be 
\[
U(\g) =\kk \langle e,\, f,\, h \rangle / 
\langle ef+fe-h, \ he-eh, \ hf-fh, \ e^2, \ f^2 \rangle,
\]
and the homogenization of $U(\g)$ by a central element $t$, analogous to 
that used in \cite{LS}, is a quadratic algebra~$H$, where 
\[
H = \frac{ \kk \langle e,\, f,\, h , \, t\rangle }%
{ \langle ef+fe-ht, \ he-eh, \ hf-fh, \ e^2, \ f^2 , \ 
et-te,\  ft-tf,\  ht-th \rangle}\, ,
\]
such that $H/H(t-1) \cong U(\g)$. 
In particular, $e$ and $f$ are zero divisors in $U(\g)$, so $H$~is not a 
domain, and hence not regular~(\cite{ATV2,Lev}). 
In spite of this, $\g \hookrightarrow U(\g)$ (cf.\ Remark~\ref{g.in.U(g)}).
In order to address the presence of zero divisors in $H$, we use a quadratic 
algebra~$\hat H$ that maps onto~$H$, namely 
\[
\hat H = \frac{\kk \langle e,\, f,\, h,\, t \rangle }%
{ \langle ef+fe-ht,\  he-eh,\  hf-fh,\  et-te,\  ft-tf,\  ht-th \rangle} \, ,
\] 
and $H = \hat H/\langle e^2, \, f^2 \rangle$.
The algebra $\hat H$ is defined in \cite{RC} and is proved therein to be an 
AS-regular algebra.  We carry over the $\Z_2$-grading from $\g$ to $\hat H$ 
by taking $|e|=|f|=\overline{1}$ and $|h|=|t|=\overline{0}$.
Furthermore, $\hat H$~has also the standard $\Z$-grading, with respect to which 
each of the generators $e,\, f,\, h,\, t$ has degree one.  

In \cite{RC}, the lines in $\PP^3$ that correspond to the right line modules 
over $\hat H$ are determined. That work entails identifying $\PP^3$ with 
$\PP((\kk e \oplus \kk f \oplus \kk h \oplus \kk t)^*)$. 
We note that the symmetry in the relations of $\hat H$ implies that the lines
corresponding to the right line modules also correspond to the left line 
modules. 

\begin{prop}\cite{RC}\label{LinVar}
The lines in $\PP^3$ that correspond to the $($left$)$ line modules over $\hat H$ 
are\\[-7mm] 
\begin{enumerate}
\item[$($a$)$] all lines in $\PP^3$ that meet the line $\V(h,\, t)$, and 
\item[$($b$)$] all lines on the quadric $\V(ht-2ef)$.
\end{enumerate}
If \,$\V(u,\,v)$ is such a line, where $u,\, v \in \hat H_1$ are linearly
independent, then the corresponding line module is isomorphic to the 
$\hat H$-module $\hat H/(\hat H u + \hat H v)$.
\end{prop}

\begin{lemma} \label{2Dsubssl11} \hfill
\quad\\[-7mm]
\begin{enumerate}
\item[$($a$)$] 
The $2$-dimensional subalgebras of $\g$ are the subspaces 
$\kk h \oplus \kk (\alpha e + \beta f)$ for all $(\alpha, \, \beta) \in \PP^1$.
\item[$($b$)$] 
All $2$-dimensional subalgebras of $\g$ are $\Z_2$-graded.
\item[$($c$)$] 
If $\h = \kk h \oplus \kk (\alpha e + \beta f)$, where $(\alpha, \, \beta) \in \PP^1$,
and if $\phi \in \h^{\ast}$, then $\phi$ determines a $1$-dimensional $\h$-module 
$\kk_{\phi}$ $($via $a \cdot y = \phi(a)y$ for all $a \in \h$, $y \in \kk_\phi$$)$
if and only if $(\phi(\alpha e + \beta f))^2 = \alpha\beta\phi(h)$. 
\end{enumerate}
\end{lemma}
\begin{pf}\hfill\\
\indent (a) For $(\alpha,\, \beta) \in \PP^1$, it is straightforward to check that 
$\kk h \oplus \kk (\alpha e + \beta f)$ is a subalgebra of $\g$ of dimension two. 

In order to show that all 2-dimensional subalgebras are of this form, let $\h$ be 
such a subalgebra, and suppose $\alpha e + \beta f + \gamma h \in \h$ for some 
$\alpha, \, \beta, \, \gamma \in \kk$. Since $\h$ is a subalgebra,  we have that 
$$
2 \alpha \beta h = 
[\alpha e + \beta f + \gamma h, \, \alpha e + \beta f + \gamma h] \in \h.
$$
Thus, if $\alpha \beta \neq 0$, then $h \in \h$, and so 
$\alpha e + \beta f \in \h$, proving that $\h$ has the desired form.  
On the other hand, if $\alpha \beta =0$ for {\it all} elements 
$\alpha e + \beta f + \gamma h$ of $\h$, then it follows that either 
$\h = \kk h \oplus \kk e$ or $\h = \kk h \oplus \kk f$. 

(b) This follows from~(a), since $\alpha e + \beta f$ is $\Z_2$-homogeneous for 
all $(\alpha , \, \beta) \in \PP^1$.

(c) Let $\phi \in \h^{\ast}$ and write $\phi(h) = \lambda \in \kk$ and 
$\phi(\alpha e + \beta f) = \gamma \in \kk$. 
The map $\phi$ determines a representation of $\h$ on $\kk$ if and only if 
\[
[x,\, y] \cdot 1 = \phi(x)\phi(y) - (-1)^{|x||y|}\phi(y)\phi(x)
\] 
for all homogeneous $x,\, y \in \h$. Checking this criterion, we obtain 
\begin{align*}
[h,\, h] \cdot 1 &= 0\cdot 1 = 0 = \lambda^2 - (-1)^0 \lambda^2,\\
[h, \,\alpha e + \beta f]\cdot 1&= 0 \cdot 1 = 0 = 
\lambda\gamma - (-1)^0 \gamma\lambda, 
\end{align*} 
whereas
\begin{align*}
[\alpha e + \beta f, \, \alpha e + \beta f] \cdot 1 &= 2 \alpha\beta h \cdot 1 
= 2\alpha\beta\lambda,
\end{align*}
and
\begin{align*}
(\phi(\alpha e + \beta f))^2 - (-1)^1
(\phi(\alpha e + \beta f))^2 &= 2 \gamma^2,
\end{align*}
which completes the proof.
\end{pf}

For a subalgebra $\h= \kk h \oplus \kk (\alpha e + \beta f)$, where 
$(\alpha, \, \beta) \in \PP^1$, and $\phi \in \h^*$, we consider the 
$\hat H$-module 
\[
L_{(\h, \, \phi)}=
\frac{\hat H}{%
\hat H (h- \phi(h) t)+ \hat H (\alpha e + \beta f - \phi(\alpha e + \beta f)t)}.
\]

\begin{thm}\label{Lhphi}
The line modules over $\hat H$ that correspond to the lines in $\PP^3 \setminus 
\V(t)$ that meet $\V(h,\, t)$ are precisely the $\hat H$-modules  
$L_{(\h, \, \phi)}$ for all $2$-dimensional subalgebras $\h$ of $\g$ and linear 
functionals $\phi \in \h^{\ast}$.
\end{thm}
\begin{pf}
If $\ell$ is a line in $\PP^3 \setminus \V(t)$ that meets $\V(h , \, t)$, 
then $\ell = \V(h - \lambda t, \, \alpha e + \beta f - \gamma t)$ for some 
$\lambda, \, \gamma \in \kk$ and $(\alpha, \, \beta) \in \PP^1$. 
By Proposition~\ref{LinVar}, such a line corresponds to a line module~$L$ 
over~$\hat H$ and
$L \cong \hat H/(\hat H(h - \lambda t) + \hat H(\alpha e + \beta f - \gamma t))$. 
To the module $L$, we associate the pair $(\h, \, \phi)$ where 
$\h = \kk h \oplus \kk(\alpha e + \beta f)$ is a $2$-dimensional subalgebra of~$\g$,
by Lemma~\ref{2Dsubssl11}(a),
and $\phi \in \h^{\ast}$ satisfies $\phi(h) = \lambda$ and 
$\phi(\alpha e + \beta f) = \gamma$. Conversely, given $(\h, \, \phi)$ as in 
the statement, and noting Lemma~\ref{2Dsubssl11}(a), the reverse process yields 
the line 
$\ell' = \V(h - \phi(h) t, \, \alpha e + \beta f - \phi(\alpha e + \beta f) t)$, 
where $(\alpha, \, \beta) \in \PP^1$ and $\h = \kk h \oplus \kk(\alpha e + \beta f)$,
and all such lines meet $\V(h,\, t)$ but do not lie on~$\V(t)$. 
By Proposition~\ref{LinVar}, 
the line module corresponding to $\ell'$ is the module $L_{(\h, \, \phi)}$.
\end{pf}

\begin{defn}\label{phi}
Let $\h = \kk h \oplus \kk (\alpha e + \beta f)$ for some $(\alpha, \, \beta) 
\in \PP^1$, and let $\phi \in \h^*$. We call $\phi$ $\Z_2$-graded if 
$\phi(\alpha e + \beta f) =  0$. 
\end{defn}

\begin{prop}\label{mainthm0}
For any $2$-dimensional subalgebra $\h$ of $\g$ and $\Z_2$-graded linear functional 
$\phi \in \h^*$, the $\hat H$-module $L_{(\h, \, \phi)}$ is a $\Z_2$-graded line 
module over $\hat H$.  Conversely, for any $\Z_2$-graded line module~$L$ over 
$\hat H$ on which $t$~acts without torsion, $L \cong L_{(\h,\, \phi)}$ for some 
$2$-dimensional subalgebra $\h$ of $\g$ and $\Z_2$-graded linear functional 
$\phi \in \h^*$. 
\end{prop}
\begin{pf}
Let $(\h, \, \phi)$ be as in Definition~\ref{phi}.  By Theorem~\ref{Lhphi},
$L_{(\h,\, \phi)}$ is a line module over~$\hat H$.  Since $\phi$ is $\Z_2$-graded,
$\alpha e + \beta f - \phi(\alpha e + \beta f)t = \alpha e + \beta f$, which is 
$\Z_2$-homogeneous.  As $h - \phi(h)t$ is also $\Z_2$-homogeneous, it follows that 
$L_{(\h,\, \phi)}$ is $\Z_2$-graded.

For the converse, suppose $L$ is a $\Z_2$-graded line module over $\hat H$ on
which $t$ acts without torsion.  By Proposition~\ref{LinVar}, 
$L \cong \hat H/\hat H \h'$ where $\h'$ is a $2$-dimensional subspace
of $\hat H_1$.  As $L$ is $\Z_2$-graded, it follows that either 
$\h' = \kk e \oplus  \kk f$ or 
$\h' = \kk h \oplus  \kk t$ or 
$\h' = \kk (\delta_1 h - \delta_2 t) \oplus  \kk (\alpha e + \beta f)$,
where $(\delta_1, \, \delta_2) , \, (\alpha, \, \beta) \in \PP^1$. 
In the first case, $L$ corresponds to the line $\V(e,\, f)$, which does not meet 
$\V(h,\, t)$ and does not lie on $\V(2ef - ht)$; hence, by Proposition~\ref{LinVar}, 
$\h' \neq \kk e \oplus \kk f$. In the second case, $t$ acts by torsion on $L$, which 
contradicts the hypothesis. It follows that 
\[
L \cong 
\frac{\hat H}{\hat H(\delta_1 h - \delta_2 t) +  \hat H(\alpha e + \beta f)},
\] 
where $(\delta_1, \, \delta_2) , \, (\alpha, \, \beta) \in \PP^1$. However, the case
that $\delta_1 = 0$ yields a module on which $t$ acts with torsion, and so we may
assume $\delta_1 = 1$.  Noting Lemma~\ref{2Dsubssl11}
and defining $\h = \kk h \oplus \kk (\alpha e + \beta f)$ and 
$\phi \in \h^*$ by $\phi (h) = \delta_2$ and  $\phi (\alpha e + \beta f) = 0$, 
completes the proof.
\end{pf}

\medskip

As mentioned earlier, our goal is to extend the results of \cite{LS,LV} to 
$\g$, where $\g = \fslll$.  
Extending the methods of \cite{LV} to our setting entails using the module $U(\g)
\otimes_{U(\h)} \kk_{\phi}$ where $\h$~is a $2$-dimensional subalgebra of~$\g$,
$U(\g)$ and $U(\h)$ are the universal enveloping algebras of $\g$ and $\h$ 
respectively, and $\kk_{\phi}$ is a 1-dimensional $\h$-module (cf.\   
Lemma~\ref{2Dsubssl11}(c)).  
Applying the methods of \cite[Section~2]{LV} in our setting produces a graded 
$H$-module from the $U(\g)$-module $U(\g) \otimes_{U(\h)} \kk_{\phi}$  using the 
central generator~$t$. The idea is as follows. One observes that $U(\g)$ is a 
filtered algebra by taking 
$U(\g)_0 = \kk$, $U(\g)_1 = \kk \oplus \g$, 
$U(\g)_2 = U(\g)_1 + \g^2$, $\ldots$,
for the filtration, and hence that 
$M = U(\g) \otimes_{U(\h)} \kk_{\phi}$  is a filtered $U(\g)$-module, where
the filtration is given by $M_i = U(\g)_i \otimes \kk_{\phi}$ for all $i$.
The graded $H$-module, $\mathcal H(M)$, obtained from this data is given by 
$\mathcal H(M) = 
\bigoplus_{i=0}^\infty M_i t^i$. In this case, $t \in H$ acts on the graded
module via \,$t\cdot((a\otimes b)t^i) = (a\otimes b)t^{i+1}$ and 
$x \in \kk e \oplus \kk f \oplus \kk h \subset H$ acts on the module via 
$x\cdot((a\otimes b)t^i) = (xa\otimes b)t^{i+1}$  for all $i$.

However, in $U(\g)$, we have that $e^2 = 0 = f^2$, and so $U(\g)$ does not have a
PBW basis that contains $e^i f^j h^k$ for \textit{all} $i,\, j,\, k \geq 0$. Thus,  
there is no reason to expect the graded $H$-module obtained via this method 
to be a line module (indeed, in general, it will not have the desired 
Hilbert series). Instead, we replace $U(\g)$, $U(\h)$ and $H$ in the above 
construction with algebras $\widehat{U(\g)}$, $\widehat{U(\h)}$, 
and $\hat H$, respectively, where $\hat H$ is given above and 
\[
\widehat{U(\g)}  = \kk \langle e,\, f,\, h \rangle / 
                      \langle ef+fe-h, \ he-eh, \ hf-fh\rangle , 
\]
and $\widehat{U(\h)}$ is the subalgebra of $\widehat{U(\g)}$ generated by $\h$.
By construction, $\widehat{U(\g)} \overset{\chi}{\onto}\ U(\g)$ (where $\chi$ is
the canonical map) and $\chi|_{\widehat{U(\h)}} : \widehat{U(\h)} \onto\ U(\h)$.
Moreover, we replace the $U(\h)$-module $\kk_{\phi}$ by a 
$\widehat{U(\h)}$-module $\widehat{\kk_{\phi}}$ where 
$\widehat{\kk_{\phi}} = \kk_{\phi}$ as a vector space, and the action is 
given by $y\cdot a = \chi(y) a$ for all 
$y \in \widehat{U(\h)}$, $a \in \widehat{\kk_{\phi}}$.  The algebra 
$\widehat{U(\g)}$ has a filtration analogous to that for $U(\g)$; namely, 
$\widehat{U(\g)}_0 = \kk$, 
$\widehat{U(\g)}_i = \widehat{U(\g)}_{i-1} + \g^i$, for all $i\geq 1$.
With these replacements, we obtain the filtered $\widehat{U(\g)}$-module 
$N = \widehat{U(\g)} \otimes_{\widehat{U(\h)}} \widehat{\kk_{\phi}}$ and the graded
$\hat H$-module $\mathcal H(N) = \bigoplus_{i=0}^\infty N_i t^i$, where 
$N_i = \widehat{U(\g)}_i \otimes_{\widehat{U(\h)}} \widehat{\kk_{\phi}}$.

\begin{rmk}\label{g.in.U(g)}
Before continuing we remark that, although $U(\g)$ and $\widehat{U(\g)}$ are not
graded algebras, the above discussion implicitly assumes that 
$\g \hookrightarrow U(\g)$  and $\g \hookrightarrow \widehat{U(\g)}$.  
These facts both follow from the observation that 
$\widehat{U(\g)}$ is the Ore extension 
$\widehat{U(\g)} = \kk[h, \, e][f; \sigma, \, \delta]$, 
where $\kk[h, \, e]$ is the polynomial ring on two generators, 
$\sigma (e) = -e$, $\sigma (h) = h$, $\delta (e) = h$, $\delta (h) = 0$, 
and so $\widehat{U(\g)}$ has PBW basis 
$\{ e^i f^j h^k : i,\, j,\, k \geq 0\}$, 
and that $U(\g) = \widehat{U(\g)}/\langle e^2, \, f^2\rangle$, a factor of 
$\widehat{U(\g)}$ by monomial relations.  These embeddings allow 
$\mathcal H(M)$ and $\mathcal H(N)$ to be well defined.
\end{rmk}


In attempting to adapt the methods of~\cite{LV} from the setting of 
$\fsl_2(\C)$ (where  $f_{\lambda}([\b, \, \b] ) = 0$  as in Remark~\ref{borel}) 
to our setting, we find that assuming the analogous condition regarding $\phi \in
\h^*$ (that is, $\phi$ is $\Z_2$-graded) typically does not yield a 
1-dimensional $\h$-module (cf.\ Lemma~\ref{2Dsubssl11}(c)).
Hence, ``forgetting'' the $\Z_2$-grading might be beneficial,
and this is the approach used in the next result.

\begin{thm}\label{mainthm1}\quad\\[-6mm]
\begin{enumerate}
\item[$($a$)$] Let $\h = \kk h \oplus \kk(\alpha e + \beta f)$ be a subalgebra of $\g$ 
where $(\alpha, \, \beta) \in \PP^1$. If $\phi \in \h^{\ast}$, where 
$(\phi(\alpha e + \beta f))^2 = \alpha\beta \phi(h)$, then 
$\mathcal H (\, \widehat{U(\g)} \otimes_{\widehat{U(\h)}} \widehat{\kk_{\phi}}\, )$
is isomorphic to the line module $L_{(\h, \, \phi)}$.
\item[$($b$)$] Let $L$ be the line module 
$L = \hat H/(\hat H(h - \lambda t) + \hat H(\alpha e + \beta f - \gamma t))$,
where $\lambda, \, \gamma \in \kk$ and  $(\alpha, \, \beta) \in \PP^1$.  If 
$\gamma^2 = \alpha\beta\lambda$, then $L \cong
\mathcal H (\, \widehat{U(\g)} \otimes_{\widehat{U(\h)}} \widehat{\kk_{\phi}}\, )$, 
where $\h = \kk h \oplus \kk(\alpha e + \beta f)$, $\phi(h) = \lambda$ and 
$\phi(\alpha e + \beta f) = \gamma$.  
\end{enumerate}
\end{thm}
\begin{pf}\hfill\\
\indent 
(a) Let $\psi$ denote the surjective map $\psi : \widehat{U(\g)} \to  
\widehat{U(\g)} \otimes_{\widehat{U(\h)}} \widehat{\kk_{\phi}}$, where
$\psi (x) = x \otimes 1$ for all $x \in \widehat{U(\g)}$.  
Since $\widehat{U(\g)}$ has a PBW basis (cf.\ Remark~\ref{g.in.U(g)}) and since 
$\widehat{U(\g)} \otimes_{\widehat{U(\h)}} \widehat{\kk_{\phi}}$ is a filtered 
$\widehat{U(\g)}$-module, we obtain $\dim_\kk(\psi(\widehat{U(\g)}_i)) = i+1$,
for all $i \geq 0$, analogous to the situation in the proof of 
\cite[Theorem~2.2(2)]{LV}.  It follows that each nonzero 
homogeneous component of degree~$i$ of
$\mathcal H (\, \widehat{U(\g)} \otimes_{\widehat{U(\h)}} \widehat{\kk_{\phi}}\, )$
has dimension $i+1$, so that the Hilbert series of 
$\mathcal H (\, \widehat{U(\g)} \otimes_{\widehat{U(\h)}} \widehat{\kk_{\phi}}\, )$
equals that of the polynomial ring on two variables. Moreover, 
$\mathcal H (\, \widehat{U(\g)} \otimes_{\widehat{U(\h)}} \widehat{\kk_{\phi}}\, )
= \hat H(1 \otimes 1)$, and so is cyclic.  Hence, 
$\mathcal H (\, \widehat{U(\g)} \otimes_{\widehat{U(\h)}} \widehat{\kk_{\phi}}\,)$
is a line module. The left annihilator in $\hat H$
of $1 \otimes 1$ contains $h - \phi(h) t$ and 
$\alpha e + \beta f - \phi (\alpha e + \beta f) t$. Hence, 
$\mathcal H (\, \widehat{U(\g)} \otimes_{\widehat{U(\h)}} \widehat{\kk_{\phi}}\, )$
is a homomorphic image of the $\hat H$-module $L_{(\h, \, \phi)}$.
By Theorem~\ref{Lhphi}, $L_{(\h, \, \phi)}$ is a line module, so 
the two modules have the same Hilbert series, and thus are isomorphic.

(b) By Proposition~\ref{LinVar}, $L$ is the line module corresponding to the 
line $\V(h - \lambda t , \, \alpha e + \beta f - \gamma t)$.
If we take $\h = \kk h \oplus \kk(\alpha e + \beta f)$ and define $\phi \in 
\h^{\ast}$ by $\phi(h) = \lambda$, $\phi(\alpha e + \beta f) = \gamma$, then,
by Lemma~\ref{2Dsubssl11}(a), $\h$  is a subalgebra of $\g$, and, 
by Lemma~\ref{2Dsubssl11}(c), the condition 
$\gamma^2 = \alpha\beta\lambda$ guarantees the existence of the $1$-dimensional 
$\h$-module~$\kk_{\phi}$ and its counterpart $\widehat{\kk_{\phi}}$,
so the result follows from~(a).
\end{pf}

In contrast with Theorem~\ref{thmLS}, where each line module on which $t$~acts 
without torsion corresponds to a line  $l_{(\b , \, \lambda)}$, and hence to a
1-dimensional $\b$-module (cf.\ Remark~\ref{borel}), the last result
suggests that, in the setting of $\fslll$,  some of the line modules 
$L_{(\h, \, \phi)}$ appear not to correspond to any $1$-dimensional $\h$-module.
This is likely a consequence of $\hat H$~``forgetting'' that $e^2 = 0 = f^2$ 
in $U(\g)$; that is, $\h$ (and hence, $U(\h)$) has fewer 1-dimensional modules 
than $\widehat{U(\h)}$~has.

\bigskip
\bigskip


\section{The Color Lie Algebra $\slc(\kk)$} \label{sec:slc}

Owing to the results used from \cite{RC}, we assume that char$(\kk) = 0$ in this
section. Our goal is to extend the results of \cite{LS,LV} to the setting of 
the color Lie algebra $\slc(\kk)$ mentioned in the Introduction.  We prove results 
for $\slc(\kk)$ that are analogous to those proved in the previous section for 
$\fslll$, with the main result of this section being Theorem~\ref{mainthm2}.

In this section, we consider the color Lie algebra $\g = \slc(\kk)$ that is 
derived from the Lie algebra $\fsl_2(\kk)$ via a process that is described 
in~\cite{CSV}. Recall that $\g$ has basis $\{a_1, \, a_2, \, a_3\}$ and 
color-Lie bracket defined by 
\[
\langle a_1, \, a_2 \rangle = a_3
= \langle a_2, \, a_1 \rangle, \quad 
\langle a_2, \, a_3 \rangle = a_1
= \langle a_3, \, a_2 \rangle, \quad 
\langle a_3, \, a_1 \rangle =a_2
= \langle a_1, \, a_3 \rangle, \quad 
\langle a_i, \, a_i \rangle =0
,
\]
for all $i$. Moreover, $\g$~is $G$-graded, where $G=\Z_2 \times \Z_2$, 
$\g = \bigoplus_{\alpha \in G} \g_{\alpha}$ and  
$$
\g_{(0,\,0)}=\{0\}, \quad \g_{(1,\,0)}=\kk a_1, \quad \g_{(0,\,1)} = \kk a_2, 
\quad \g_{(1,\,1)} = \kk a_3.
$$
As in~\cite{CSV}, $\g$ can be viewed as a cocycle twist of~$\fsl_2(\kk)$. 
Continuing our terminology from the previous section, we refer to any subspace 
of~$\g$ that is closed under the color-Lie bracket as a {\em subalgebra} of~$\g$.

The universal enveloping algebra, $U(\g)$, of~$\g$ is defined to be 
\[
U(\g) = \frac{\kk \langle a_1, \, a_2, \, a_3 \rangle}%
{ \langle a_1 a_2 + a_2 a_1 - a_3, \  a_2 a_3 + a_3 a_2 - a_1, \  
a_3 a_1 + a_1 a_3 - a_2 \rangle }.
\]
We denote the homogenization of $U(\g)$ by a single central element~$a_4$ by 
the algebra 
\begin{align*}
H &= \kk \langle a_1, \, a_2, \, a_3, \, a_4 \rangle / 
\langle 
a_1 a_2 + a_2 a_1 - a_3 a_4, \ 
a_2 a_3 + a_3 a_2 - a_1 a_4, \\
&\hspace{.35in}
a_3 a_1 + a_1 a_3 - a_2 a_4, \  
a_1 a_4 -a_4 a_1, \ 
a_2 a_4 - a_4 a_2, \ 
a_3 a_4 - a_4 a_3 \rangle.
\end{align*}
If we define $|a_4|=(0,\, 0) \in G$, then $H$ is a $G$-graded algebra.
However, Theorem~\ref{mainthm2} below suggests that the $G$-grading contributes
nothing towards accomplishing our objective.

We may recover $U(\g)$ from $H$ via $U(\g) \cong H/H(a_4-1)$. Using Bergman's Diamond 
Lemma, it is straightforward to see that $U(\g)$ has PBW basis 
$\{ a_1^i a_2^j a_3^k : i,\, j,\, k \geq 0\}$; thus, $\g
\hookrightarrow U(\g)$, and $\h \hookrightarrow U(\h) \hookrightarrow U(\g)$
for all subalgebras $\h$ of $\g$.

We identify $\PP^3$ with 
$\PP((\kk a_1 \oplus \cdots \oplus \kk a_4)^*)$. 
Owing to the symmetry of the relations defining~$H$, the right line modules 
over~$H$ and the left line modules over~$H$ are parametrized by the same lines
in $\PP^3$.  By \cite{RC}, $H$ is AS-regular and the line modules over $H$ are 
given by the following result.

\begin{prop}\cite{RC}\label{lineSchemesl2c} 
The $($left$)$ line modules over~$H$ are given by the following $13$ types of 
lines in~$\PP^3$, and conversely:
\begin{enumerate}
\item[$1(a)$] 
all lines in $\V(a_1 + a_2)$ that pass through $(1,\, -1,\, 0,\, 0)$,
\item[$1(b)$] 
all lines in $\V(a_1 - a_2)$ that pass through $(1, \, 1, \, 0, \, 0)$,
\item[$2(a)$] 
all lines in $\V(a_1 + a_3)$ that pass through $(1, \, 0, \, -1, \, 0)$,
\item[$2(b)$] 
all lines in $\V(a_1 - a_3)$ that pass through $(1, \, 0, \, 1, \, 0)$,
\item[$3(a)$] 
all lines in $\V(a_2 + a_3)$ that pass through $(0, \, 1, \, -1, \, 0)$,
\item[$3(b)$] 
all lines in $\V(a_2 - a_3)$ that pass through $(0, \, 1, \, 1, \, 0)$,
\item[$4(a)$] 
all lines in $\V(a_4 - 2 a_3)$ that pass through $(1,\, -1,\, 0,\, 0)$,
\item[$4(b)$] 
all lines in $\V(a_4 + 2 a_3)$ that pass through $(1, \, 1, \, 0, \, 0)$,
\item[$5(a)$] 
all lines in $\V(a_4 - 2 a_2)$ that pass through $(1, \, 0, \, -1, \, 0)$,
\item[$5(b)$] 
all lines in $\V(a_4 + 2 a_2)$ that pass through $(1, \, 0, \, 1, \, 0)$,
\item[$6(a)$] 
all lines in $\V(a_4 - 2 a_1)$ that pass through $(0, \, 1, \, -1, \, 0)$,
\item[$6(b)$] 
all lines in $\V(a_4 + 2 a_1)$ that pass through $(0, \, 1, \, 1, \, 0)$,
\item[$7$.] all lines in $\V(a_4)$.
\end{enumerate}
\end{prop}

As in Section~\ref{sec:sl11}, if $\h$ is a subalgebra of $\g$, then we denote an
$\h$-module of dimension one by~$\kk_\phi$, where $\phi \in \h^*$  and 
$a \cdot y = \phi(a)y$ for all $a \in \h$, $y \in \kk_\phi$. The next result shows
that not all elements of $\h^*$ give rise to an $\h$-module.


\begin{lemma}\label{2dimlColorLie} \hfill
\begin{enumerate}
\item[$(a)$] 
The $2$-dimensional subalgebras of $\g$ are precisely the subspaces 
$\kk a_i \oplus \kk(a_j + \mu a_k)$ for all distinct $i,\, j,\, k$ 
and $\mu = \pm 1 \in \kk$.
In particular,  no $2$-dimensional subalgebra of $\g$ is $G$-graded.
\item[$(b)$] 
Let $\mu = \pm 1 \in \kk$ and $\{ i, \, j, \, k\} = \{ 1, \, 2, \, 3\}$.
If $\h$ is the subalgebra $\h = \kk a_i \oplus \kk(a_j + \mu a_k)$, then 
$\h$ has exactly two one-parameter families of \,$1$-dimensional 
modules$:$
\[ 
\{\kk_\phi :  \phi \in \h^*, \ \phi (a_j + \mu a_k) = 0 \} 
\qquad \text{and} \qquad 
\{\kk_\phi :  \phi \in \h^*, \ \phi (a_i) = \mu/2 \}.
\]
\end{enumerate}
\end{lemma}
\begin{pf}\hfill\\
\indent (a) 
We first observe that any subspace of the given form is a subalgebra. 
Conversely, suppose $\h$ is a $2$-dimensional subalgebra and write 
$\h = \kk(\alpha_1 a_1 + \alpha_2 a_2 + \alpha_3 a_3) \oplus 
\kk(\beta_1 a_1 + \beta_2 a_2 + \beta_3 a_3)$, 
where $(\alpha_1, \, \alpha_2, \, \alpha_3)$,  
$(\beta_1, \, \beta_2, \, \beta_3) \in \PP^2$. 
Since $\dim(\h) = 2$, the matrix $X$ of coefficients, where  
\[ 
X = \begin{bmatrix}
\alpha_1 & \alpha_2 & \alpha_3\\[1mm]
\beta_1 & \beta_2 & \beta_3
\end{bmatrix},
\] 
has rank two.  By symmetry and row operations on $X$, it follows that we may assume 
$\h = \kk v_1 \oplus \kk v_2$, where $v_1 = a_1 + \alpha a_3$ and 
$v_2 = a_2 + \beta a_3$,  for some $\alpha, \, \beta \in \kk$. 
Since $\h$ is a subalgebra, 
%
\[ 
2\alpha a_2 = \langle v_1, \, v_1 \rangle \in \h,\quad
2\beta a_1 = \langle v_2, \, v_2 \rangle \in \h\quad
\text{and} \quad 
\alpha a_1 + \beta a_2 + a_3 = \langle v_1, \, v_2 \rangle \in \h.
\]
However, $\dim(\h) = 2$, so it follows that the matrix of coefficients of the 
elements $v_1$, $v_2$, $2\alpha a_2$, $2\beta a_1$, $\alpha a_1 + \beta a_2 + a_3$ 
has rank at most two; that is, all $3 \times 3$ minors of the matrix 
\[
\begin{bmatrix}
\ 1 & 0 & 0 & 2\beta & \alpha\ \\[1mm]
\ 0 & 1 & 2\alpha & 0 & \beta\ \\[1mm]
\ \alpha & \beta & 0 & 0 & 1\  
\end{bmatrix} 
\] 
are zero.  Hence, $2\alpha\beta = 0 = \alpha^2 + \beta^2 - 1$, which implies
$(\alpha, \, \beta) = (0, \, \pm 1)$ or $(\pm 1, \, 0)$ and so $\h$ is of the 
desired form. For $\mu = \pm 1$, $a_j + \mu a_k$ is not $G$-homogeneous
for all distinct $j$ and $k$, so $\h$ is not $G$-graded.

(b) By symmetry, it suffices to consider $(i,\, j,\, k) = (3,\, 1,\, 2)$.
Let $\phi \in \h^*$.
There is only one relation that needs to be checked, namely 
\[
a_3\cdot ((a_1 + \mu a_2)\cdot 1) + 
(a_1 + \mu a_2)\cdot (a_3 \cdot 1) 
= \mu (a_1 + \mu a_2)\cdot 1.
\]
Writing $\phi(a_3) = \alpha \in \C$ and $\phi(a_1 + \mu a_2) = \beta \in \C$, and 
evaluating each side of the relation implies that 
$\alpha \beta + \beta \alpha = \mu \beta$; that is, 
$\kk$ is an $\h$-module if and only if $\beta (2 \alpha - \mu) = 0$. This latter
equation has solution set
$\{ (\alpha' , \, 0), \  (\mu/2, \, \beta') : \alpha', \, \beta' \in \kk\}$,
so the result follows.
\end{pf}
%

Given $\h$ and $\kk_\phi$ as in Lemma~\ref{2dimlColorLie}, we may consider 
the $U(\g)$-module $U(\g) \otimes_{U(\h)} \kk_\phi$ and the $H$-module 
$\mathcal H (\, U(\g) \otimes_{U(\h)} \kk_\phi \, )$, where the latter is 
defined in a manner analogous to that described in \cite[Section~2]{LV}  and 
in Section~\ref{sec:sl11} of this article. Moreover, since $U(\g)$ has a 
PBW basis, arguments similar to those provided in Section~\ref{sec:sl11} prove 
the next result.


\begin{thm}\label{mainthm2}  
For each $k \in \{ 1, \, 2, \, 3 \}$, the $a_k$-torsion $($left$)$ line modules 
over~$H$, on which $a_4$~acts without torsion, are homogenizations of 
induced modules. More precisely, we have the following.
\begin{enumerate}
\item[$(a)$] 
There exists a one-to-one correspondence between all pairs $(\h, \, \phi)$, such 
that  $\h =  \kk a_i \oplus \kk(a_j + \mu a_k)$, $\{  i, \, j, \, k \} = 
\{1,\, 2,\,3\}$, $\mu = \pm 1 \in \kk$, and $\phi \in \h^*$ where 
$\phi (a_j + \mu a_k) = 0$, and line modules given by 1$(a)$-3$(b)$ in
Proposition~\ref{lineSchemesl2c} on which $a_4$ acts without torsion$;$ this
correspondence is given by  
\[
\mathcal H (\, U(\g) \otimes_{U(\h)} \kk_\phi \, ) \ \cong \  
\frac{H}{H(a_i - \phi(a_i) a_4) + H (a_j + \mu a_k)}\,.
\]
\item[$(b)$] 
There exists a one-to-one correspondence between all pairs $(\h, \, \phi)$, such 
that  $\h = \kk a_i \oplus \kk(a_j + \mu a_k)$, $\{  i, \, j, \, k \} = 
\{1,\, 2,\,3\}$, $\mu = \pm 1 \in \kk$, and $\phi \in \h^*$ where
$\phi (a_i)  = \mu/2$, and line modules given by 4$(a)$-6$(b)$ in
Proposition~\ref{lineSchemesl2c} on which $a_4$ acts without torsion$;$ this
correspondence is given by  
\[
\mathcal H (\, U(\g) \otimes_{U(\h)} \kk_\phi \, ) \ \cong \  
\frac{H}{H(a_4 - 2 \mu a_i) + H (a_j + \mu a_k - \phi(a_j + \mu a_k) a_4)}\,.
\]
\quad\\[-13mm]
\qedsymbol
\end{enumerate}
\end{thm}

This last result suggests that the $G$-grading on~$\g$ and the 
$G$-grading on~$H$ play no role in dictating the correspondence between 
line modules over $H$ and pairs $(\h, \, \phi)$ as discussed in the statement.



\section{The Lie Superalgebra $\fsl(2|1)$}

In this section, we consider the Lie superalgebra $\fsl(2|1)$, 
and return to the assumption that char$(\kk) \neq 2$.
We show that, in this setting, a simple analogue of the algebra $\hat H$
(respectively, $H$) that was used in Section~\ref{sec:sl11} 
(respectively, Section~\ref{sec:slc}) yields an algebra with zero divisors and so
is not a regular algebra.

Recall that $\fsl(2|1)$ consists of the matrices $(\alpha_{ij}) \in M(3,\, \kk)$ 
with supertrace $\alpha_{11} + \alpha_{22} - \alpha_{33} = 0$,  where the 
Lie (super)bracket is defined as in Section~\ref{sec:sl11}, and the 
$\Z_2$-grading on $\fsl(2|1)$ may be described as follows. 
Using the elementary matrices, $E_{11}, E_{12}, \dots, E_{33}$, as a 
basis for $M(3,\, \kk)$, let $x_1, \ldots, x_4, y_1, \ldots, y_4$ be given by 
$x_1 = E_{11} + E_{33}$, $x_2 = E_{22} + E_{33}$, $x_3 = E_{12}$, $x_4 = E_{21}$, 
$y_1 = E_{13}$, $y_2 = E_{31}$, $y_3 = E_{23}$, $y_4 = E_{32}$.  With this
notation, $\fsl (2|1) = (\fsl (2|1))_{\overline 0} \oplus 
(\fsl (2|1))_{\overline 1}$, where 
$(\fsl (2|1))_{\overline 0} = \bigoplus_{i=1}^4 \kk x_i$ and 
$(\fsl (2|1))_{\overline 1} = \bigoplus_{i=1}^4 \kk y_i$. 

Homogenizing the universal enveloping algebra of $\fsl(2|1)$ by using a central 
element $t$ and then deleting the relations $y_i^2 = 0$, for all $i = 1, \ldots,
4$, yields a quadratic algebra $\hat H$ on the nine generators 
$x_1, \ldots, x_4, y_1, \ldots, y_4, \, t$ with 
$36$ (i.e., 
{\tiny $\begin{pmatrix}9\\[-2pt]2\end{pmatrix}$})
defining relations, among which are the relations 
$x_3 y_1 = y_1 x_3$, $y_1 y_3 = - y_3 y_1$ and $x_3 y_3 - y_3 x_3 = y_1 t$. 
These three defining relations imply that
\begin{align*}
y_1^2 t - x_3y_1y_3 
&= y_1^2t - y_1x_3y_3\\
&= y_1(y_1t - x_3y_3)\\
&= -y_1y_3x_3\\
&= y_3y_1x_3\\
&= y_3x_3y_1\\
&= (x_3y_3 - y_1t)y_1\\
&= x_3y_3y_1 - y_1^2t\\
&= -y_1^2t - x_3y_1y_3,
\end{align*}    
from which it follows that $2 y_1^2 t = 0$ in $\hat H$. Thus, $\hat H$ is not a 
domain, and hence not regular~(\cite{ATV2,Lev}). 

Consequently, if results analogous to those in Sections~\ref{sec:sl11} 
and \ref{sec:slc} hold for $\fsl(2|1)$, then either an alternative graded algebra
will need to be used, or a careful analysis with the zero divisors will need to be
performed in order to discuss the line modules.

\bigskip
\bigskip



\bigskip
\bigskip
\bigskip

{\footnotesize
\noindent
\textsc{%
S.J.~Sierra, School of Mathematics, University of Edinburgh, Edinburgh EH9 3FD, 
UK.}\\
\textit{e-mail address}:
\texttt{s.sierra@ed.ac.uk}
\\[-1mm]

\noindent
\textsc{%
\v{S}.~\v{S}penko, Departement Wiskunde, Vrije Universiteit Brussel, 
Pleinlaan 2, B-1050 Brussels, Belgium. }\qquad
\textit{e-mail address}:
\texttt{spela.spenko@vub.be}
\\[-1mm]

\noindent
\textsc{%
M.~Vancliff, Department of Mathematics, Box 19408, University of Texas at
Arlington, \\Arlington, TX 76019-0408, USA.}\qquad
\textit{e-mail address}:
\texttt{vancliff@uta.edu}
\\[-1mm]

\noindent
\textsc{%
P.~Veerapen, Department of Mathematics, Tennessee Technological University, 
Cookeville, TN 38505, USA.}\qquad
\textit{e-mail address}:
\texttt{pveerapen@tntech.edu}
\\[-1mm]

\noindent
\textsc{%
E.~Wiesner, Department of Mathematics, Ithaca College, Ithaca, NY 14850, 
USA.}\\
\textit{e-mail address}:
\texttt{ewiesner@ithaca.edu} 

}

\end{document}